\documentclass[11pt]{amsart}
\usepackage{amsmath}
\usepackage{amssymb}
\input amssym.def
\input amssym.tex

\font\tengothic=eufm10
\font\sevengothic=eufm7
\newfam\gothicfam
      \textfont\gothicfam=\tengothic
      \scriptfont\gothicfam=\sevengothic
\def\goth#1{{\fam\gothicfam #1}}

\numberwithin{equation}{section}

\newenvironment{mydef}{{\bf Definition.}}{\hspace*{\fill} \par\vspace{1ex}}

\begin{document}
\setlength{\baselineskip}{1.7em}
\newtheorem{thm}{\bf Theorem}[section]
\newtheorem{pro}[thm]{\bf Proposition}
\newtheorem{claim}[thm]{\bf Claim}
\newtheorem{lemma}[thm]{\bf Lemma}
\newtheorem{cor}{\bf Corollary}[thm]

\newcommand{\mP}{{\mathbb P}}
\newcommand{\z}{{\mathbb Z}}
\newcommand{\nn}{{\mathbb N}}
\newcommand{\R}{{\mathcal R}}
\newcommand{\A}{{\mathcal A}}
\newcommand{\x}{{\mathbb X}}
\newcommand{\y}{{\mathbb Y}}
\newcommand{\ix}{I_{\mathbb X}}
\newcommand{\mi}{{\mathcal I}}
\newcommand{\B}{{\bf B}}
\newcommand{\bi}{{\bf I}}
\newcommand{\V}{{\bf V}}
\newcommand{\cv}{{\mathcal V}}
\newcommand{\bL}{{\bf L}}
\newcommand{\calo}{{\mathcal O}}
\newcommand{\cl}{{\mathcal L}}
\newcommand{\cm}{{\mathcal M}}
\newcommand{\cn}{{\mathcal N}}
\newcommand{\m}{-\!\!--\!\!\rightarrow}
\newcommand{\smap}{\rightarrow\!\!\!\!\!\rightarrow}
\newcommand{\sfrac}[2]{\frac{\displaystyle #1}{\displaystyle #2}}
\newcommand{\gr}{\Gamma_{d+1}}
\newcommand{\img}{\Lambda_{d+1}}
\newcommand{\under}[1]{\underline{#1}}
\newcommand{\ov}[1]{\overline{#1}}
\newcommand{\sheaf}[1]{{\mathcal #1}}
\newcommand{\lex}{{\le}_{\mbox{\scriptsize lex}}}

\title{On the Rees algebra of certain codimension two perfect ideals.}
\author{H\`a Huy T\`ai}
\thanks{Current address: Institute of Mathematics, P.O. Box 631 B\`o H\^o, H\`a N\^oi 10000, Vietnam.}
\thanks{Email: tai@hanimath.ac.vn}
\address{Department of Mathematics and Statistics \\ Queen's University, Kingston ON K7L 3N6, Canada} \email{haht@mast.queensu.ca}
\maketitle

\setcounter{section}{-1}
\section{Introduction.}

Let $R$ be a commutative ring with identity and $I \subseteq R$ an ideal of $R$. The Rees algebra of $I$ is defined to be the subring 
\[ R[I t] = R \oplus It \oplus I^2t^2 \oplus \ldots ~ \subseteq R[t], \]
and denoted by $\R_R(I)$ (or simply $\R(I)$, if there is no confusion on the ring $R$ being discussed). 

The Rees algebra of an ideal is a classical object that has been studied throughout many decades. Our interest to Rees algebras comes from the fact that they provide the algebraic realizations for certain class of rational $n$-folds, namely those obtained from $\mP^n$ by blowing up at a subscheme. In this paper, we study the Rees algebras of certain codimension two perfect ideals. To be more precise, we study the Rees algebra of the defining ideal of a set of points in $\mP^2$. 

The first important result on the Rees algebras of codimension two perfect ideals in a polynomial ring was due to Morey and Ulrich (\cite[Theorem 1.3]{m-u}). They showed that the Rees algebra of  any codimension two perfect ideal with a linear presentation matrix in a polynomial ring of $d$ variables, which is minimally generated by more than $d$ generators and satisfies the condition $G_d$ (condition $G_s$, $s$ an integer, for an ideal $I$ of a ring $R$ means that the minimal number of generators of $I_p$ is less than or equal to the dimension of $R_p$ for every prime ideal $p \supseteq I$ such that $\dim R_p \le s-1$), is Cohen-Macaulay and generated by the maximal minors of a matrix of linear forms. When reduced to the class of defining ideals for a generic set of points in $\mP^2$, their restriction on the presentation matrix of the ideal requires the number of points to be a binomial coefficient number. When the number of points is arbitrary, or the set of points are not in generic position, not much is known about the Rees algebra of its defining ideal. 

In our approach to the problem, inspired by works of Mumford (\cite{mumford}) and Green (\cite{green}), we look at an arbitrary set of points in $\mP^2$ and study how the Rees algebras of the ideals generated by the homogeneous pieces of its defining ideal behave asymptotically. More precisely, suppose $\x$ is a set of points in of $\mP^2$ and $\ix = \oplus_{t \ge \alpha} I_t \subseteq R = \goth{k}[w_1,w_2,w_3]$ its defining ideal, we will study the asymptotic behaviour of the Rees algebras $\R(I_t)$ of the ideals generated by $I_t$, as $t$ gets large. These Rees algebras have a geometric significance in the following sense. Suppose $t \ge \alpha$ and $G_0, \ldots, G_m$ is a minimal system of generators for the vector space $I_t$. We consider the rational map:
\[ \varphi: \mP^2 \m \mP^m, \]
given by sending each point $P \in \mP^2 \backslash \x$ to the point $[G_0(P): \ldots: G_m(P)] \in \mP^m$. Suppose $\Gamma$ and $\ov{\Gamma}$ are the graph and its closure of $\varphi$ in $\mP^2 \times \mP^m$. It is known (cf. \cite[Chapter 2]{ha-thesis}) that the Rees algebra $\R(I_t)$ is the bi-graded coordinate ring of $\ov{\Gamma}$, and when $t$ is at least as large as the degrees of the generators of a minimal system of generators for $\ix$, $\ov{\Gamma}$ is the blowup of $\mP^2$ along the points of $\x$, embedded into the product space $\mP^2 \times \mP^m$.

The basic outline of the paper is as follows. In the first section, we start by looking at the case when the points in $\x$ are in generic position. We prove the following results. 

\begin{thm} [Theorem \ref{3-simple-2}] \label{intro-3}
Let $I = \oplus_{t \ge d} I_t$ be the defining ideal of $s = {d+1 \choose 2}$ points in
$\mP^2$ which are in generic position. Then the defining equations for the Rees algebra
$\R(I_{d+1})$ are the $2 \times 2$ minors of a $3 \times (d+2)$ matrix of linear forms.
Moreover, $\R(I_{d+1})$ is Cohen-Macaulay, and has the same Betti numbers as that of the
ideal of $2 \times 2$ minors of a generic $3 \times (d+2)$ matrix. 
\end{thm}

\begin{thm} [Theorem \ref{3-hard-2}] \label{intro-4}
Let $I = \oplus_{t \ge d} I_t$ be the defining ideal for a set of $s = {d+1 \choose 2} +
k$ ($1 \le k \le d$) points in $\mP^2$. Then for a general choice of the points, the Rees
algebra $\R(I_{d+1})$ is Cohen-Macaulay and its defining ideal is generated by the $3\times3$ minors of a
$k\times3$ matrix $B$ of linear forms, the $2\times2$ minors of a $3\times(d-k+2)$ matrix
$X$ of indeterminates and the entries of the product matrix $B.X$.  
\end{thm}

Note that these results can be extended to a larger class of codimension two
perfect ideals of any polynomial ring (Theorems \ref{3-simple-ext} and \ref{3-hard-ext}).
These results also provide the necessary information to completely answer questions on the
defining equations of certain projective embeddings of $\mP^2(\x)$, as discussed in \cite[Chapter 4]{ha-thesis}.

Theorems \ref{intro-3} and \ref{intro-4} give motivation toward the study of asymptotic
behaviour of the Rees algebras $\R(I_t)$ as $t$ gets large, for the defining ideal $\ix
= \oplus_{t \ge \alpha} I_t$ of an arbitrary set of points $\x \subseteq \mP^2$.
This study is carried out in the last section of this paper. The main result of this section is the following Mumford-typed theorem.    

\begin{thm} [Theorem \ref{3-cm}] \label{intro-1}  
Suppose $\x = \{ P_1, \ldots, P_s \}$ is an arbitrary set of $s$ points in $\mP^2$, and $\ix = \oplus_{t \ge \alpha} I_t \subseteq R = \goth{k}[w_1,w_2,w_3]$ is its defining ideal. Then, there exists an integer $d_0$ such that for all $t \ge d_0$, the Rees algebra $\R(I_t)$ of the ideal generated by $I_t$ is Cohen-Macaulay, and its defining ideal is generated by quadratics.  
\end{thm}

In this section, we also introduce the notion of being {\it arithmetic Cohen-Macaulay} (a.CM) for a subscheme of the product scheme $\mP^n \times \mP^m$. We give a characterization for this property in the following theorem.

\begin{thm}[Theorem \ref{3-cm-condition}]
Suppose $V \subseteq \mP^n \times \mP^m$ is a proper closed subscheme of dimension $d$ of $\mP^n \times \mP^m$. Then,
\renewcommand{\labelenumi}{(\arabic{enumi})}
\begin{enumerate}
\item If $V$ is a.CM, then $H^i(\mi_V(a,b)) = 0$ for all $a,b \in \z$ and $1 \le i \le d$, where $\mi_V$ is the ideal sheaf of $V$ in $\mP^n \times \mP^m$.
\item Suppose $d \not= n, m$, and $H^i(\mi_V(a,b)) = 0$ for all $a,b \in \z$ and $1 \le i \le d$. If in addition, $H^{d+1}(\mi_V(a,b)) = 0$ for all $a,b \ge 0$, and for every $j > 0$,
\[ R^j {\pi_1}_{*}(\calo_V(p,q)) = 0 \ \forall p \in \z, q \ge 0, \]
and
\[ R^j {\pi_2}_{*}(\calo_V(p,q)) = 0 \ \forall q \in \z, p \ge 0, \]
then $V$ is a.CM.
\end{enumerate}
\end{thm}

Here, $\pi_1$ and $\pi_2$ are the two projection maps $\mP^n \times \mP^m \rightarrow \mP^n$ and $\mP^n \times \mP^m \rightarrow \mP^m$ restricted to $V$.

Throughout this paper, $\goth{k}$ will be our ground field. For simplicity, we assume that $\goth{k}$ is algebraically closed and of charateristic 0. For each $n \ge 1$, we also let $\mP^n = \mP^n_{\goth{k}}$ be the $n$-dimensional projective space over $\goth{k}$.

\section{Ideal of a generic set of points} \label{3-sec-gen}

We start by considering the situation where our set of points is in generic position. In particular, let $\x = \{ P_1,\ldots,P_s \}$ be a set of $s$ distinct points in $\mP^2$ which are in generic position. Let $\ix = \oplus_{t \ge \alpha}I_t$ be the defining ideal of $\x$ in $R=\goth{k}[w_1,w_2,w_3]$, and $\mP^2(\x)$ the blowup of $\mP^2$ centered at $\x$. We now proceed by considering different cases depending on the number of points in $\x$. Even though the heart of the matter lies in the case when the number of points in $\x$ is arbitrary, we start with the simplest situation, where we have a binomial coefficient number of points.

\begin{center}
{\bf Binomial coefficient number of points}
\end{center}

Our main result  in this subsection (Theorem \ref{3-simple-2}) can be obtained from that of Morey and Ulrich (\cite[Theorem 1.3]{m-u}, our method, nevertheless, is quite different. Our argument in this case is very similar to the argument on the Room surfaces of \cite{ge-gi}. We refer the readers to Theorem 1.2 of \cite{ge-gi}. Suppose $s={d+1 \choose 2}$ (for some integer $d$). Then, from \cite{g-m}, $\sigma(\ix) = d$ and $\ix$ is generated by $I_d$. By the Hilbert-Burch theorem, these generators are the $d \times d$ minors of a $d \times (d+1)$ matrix, say $\bL$, of linear forms :
\[ \bL = ( L_{ij} ), ~ L_{ij} \in R_1 \mbox{ for } i=1,2,\ldots,d \mbox{ and } j=1,2,\ldots,d+1. \]
In this notation, 
\[ \ix = (F_1,\ldots,F_{d+1}), ~ F_i = (-1)^{i+1}\mbox{det}(\bL ~ \backslash ~ i^{\mbox{\footnotesize th}} \mbox{ column}). \]

We shall now establish the defining equations for the Rees algebra $\R(I_{d+1})$ of the ideal generated by $I_{d+1}$. 

A system of generators of the vector space $I_{d+1}$ is given by $3(d+1)$ forms $w_iF_j$ for $i=1,2,3$ and $j=1,\ldots,d+1$. Consider the rational map :
\[ \varphi_{d+1} : \mP^2 \m \mP^N, ~ N = 3(d+1) - 1, \]
given by $\varphi_{d+1}(P) = [w_iF_j]$ for any point $P \in \mP^2 \backslash \x$. Let $\gr$ and $\img$ be the graph and the image of $\varphi_{d+1}$, and let $\ov{\gr}$ and $\ov{\img}$ be their closures in appropriate spaces. We use homogeneous coordinates $[x_{ij}]_{1 \le j \le d+1, 1 \le i \le 3}$ of $\mP^N$ such that
\begin{eqnarray}
\varphi_{d+1}([\ov{w_1}:\ov{w_2}:\ov{w_3}]) & = & [\ov{x_{ij}}], \mbox{ where } \ov{x_{ij}} = \ov{w_i}F_j. \label{3-sim0}
\end{eqnarray}

The vector space dimension of $I_{d+1}$ is $2d+3$, so there must be $d$ linear dependence relations between the $w_iF_j$'s. Those relations can be found by expanding the following zero determinants:
\[ 0 = \mbox{ det } \left( 
\begin{array}{c}L_{l1} ~ L_{l2} ~ \ldots ~ L_{l, d+1} \\ \bL \end{array}
 \right) = \sum_{j=1}^{d+1} L_{lj}F_j, \]
for each $l=1,2,\ldots,d$. Now, let $L_{lj} = \sum_{i=1}^3 \lambda_{lji}w_i$, then by grouping similar terms, we get a collection of dependence relations among the $w_iF_j$'s as follows:
\[ \sum_{j=1}^{d+1} \sum_{i=1}^{3} \lambda_{lji}w_iF_j = 0, \ \forall \ l=1,2,\ldots,d. \]
This gives rise to a collection of equations, where the coordinates of the points in $\gr$ satisfy:
\begin{eqnarray}
\sum_{1 \le i \le 3, 1 \le j \le d+1} \lambda_{lji}x_{ij} & = & 0, \ \forall \ l=1,2,\ldots,d. \label{3-sim1}
\end{eqnarray}

There are exactly $d$ equations, and as it was also proved in \cite{ge-gi}, those equations are linearly independent, so they are indeed all the equations obtained from the linear dependence relations of the $w_iF_j$s.

Consider the matrix
\[ M = \left[ \begin{array}{ccccc} w_1 & x_{11} & x_{12} & \ldots & x_{1, d+1} \\
w_2 & x_{21} & x_{22} & \ldots & x_{2, d+1} \\
w_3 & x_{31} & x_{32} & \ldots & x_{3, d+1} \end{array} \right] \]

From $(\ref{3-sim0})$, it is easy to see that the points of $\gr$ satisfy all the $ 2 \times 2$ minors of $M$. Denote the collection of these equations by (**). Let 
\[ M' = \left[ \begin{array}{cccc} x_{11} & x_{12} & \ldots & x_{1, d+1} \\
x_{21} & x_{22} & \ldots & x_{2, d+1} \\
x_{31} & x_{32} & \ldots & x_{3, d+1} \end{array} \right], \]
and recall Proposition 1.1 of \cite{ge-gi}.

\begin{pro} \label{3-ge-gi}
For each $Q=[\ov{x_{ij}}] \in \mP^N$ satisfying the equations in $(\ref{3-sim1})$ and the $2 \times 2$ minors of $M'$, there exists a unique $P'=[\ov{w_1}:\ov{w_2}:\ov{w_3}] \in \mP^2$ such that the coordinates of $P'$ and $Q$ satisfy
\[ (\dag) \left\{ \begin{array}{ccc} \ov{x_{2j}} \ \ov{w_1} - \ov{x_{1j}} \ \ov{w_2} & = & 0 \\ \ov{x_{3j}}\ \ov{w_2}-\ov{x_{2j}}\ \ov{w_3} & = & 0 \\ 
\ov{x_{1j}}\ \ov{w_3} - \ov{x_{3j}}\ \ov{w_1} & = & 0 \end{array} \right. \mbox{ for } j=1,2,\ldots,d+1. \]
\end{pro}

Similar to what was done in \cite{ge-gi}, we have the following result.

\begin{thm} \label{3-simple-1}
Equations in $(\ref{3-sim1})$ and (**) are the defining equations of $\ov{\gr}$ in $\mP^2 \times \mP^N$.
\end{thm}

\begin{proof} Let $\V$ be the algebraic set in $\mP^2 \times \mP^N$ defined by all the bi-homogeneous equations in $(\ref{3-sim1})$ and (**). We shall first prove that $\V = \ov{\gr}$ as sets. 

Clearly, the coordinates of the points of $\gr$ satisfy all the equations in $(\ref{3-sim1})$ and (**), so as sets, $\gr \subseteq \V$, hence $\ov{\gr} \subseteq \V$. To prove the reverse inclusion, let $(P,Q) \in \V$ (where $P \in \mP^2$ and $Q \in \mP^N$). The coordinates of $Q$ satisfy equations in $(\ref{3-sim1})$ and the $2 \times 2$ minors of $M'$, so by Proposition \ref{3-ge-gi}, and following the same argument as that of \cite{ge-gi}, there exists a unique $P'=[\ov{w_1}:\ov{w_2}:\ov{w_3}]$ such that the coordinates of $P'$ and $Q$ satisfy $(\dag)$, and $Q$ must have the form
\begin{eqnarray}
Q = [\ov{w_1}c_1:\ov{w_2}c_1:\ov{w_3}c_1:\ldots:\ov{w_1}c_{d+1}:\ov{w_2}c_{d+1}:\ov{w_3}c_{d+1}], \label{3-break-1}
\end{eqnarray}
for some $c_1,c_2,\ldots,c_{d+1} \in \goth{k}$. This and the equations in $(\ref{3-sim1})$ implies that
\[ \bL(P') \left[ \begin{array}{c} c_1 \\ \vdots \\ c_{d+1} \end{array} \right] = \left[ \begin{array}{c} 0 \\ \vdots \\ 0 \end{array} \right]. \]
Thus, if $P' \not\in \x$ then
\begin{eqnarray}
\left[ \begin{array}{c} c_1 \\ \vdots \\ c_{d+1} \end{array} \right] = \rho \left[ \begin{array}{c} F_1(P') \\ \vdots \\ F_{d+1}(P') \end{array} \right], \label{3-break-2}
\end{eqnarray}
for some $\rho \in \goth{k}$; and otherwise, if $P' \in \x$, then $Q$ lies on the exceptional line corresponding the the blowup at $P$. Thus, $Q \in \ov{\img}$.

We also note that the equations in $(\dag)$ and the $2 \times 2$ minors of $M'$ are exactly the $2 \times 2$ minors of $M$. Thus, Proposition \ref{3-ge-gi} shows that for each $Q \in \ov{\img}$, there exists a unique $P' \in \mP^2$ such that the coordinates of $P'$ and $Q$ satisfy the $2 \times 2$ minors of $M$. This implies $P = P'$. Since the linear system $I_{d+1}$ is very ample (\cite{d-g}), the projection map $\ov{\gr} \rightarrow \ov{\img}$ defined by sending $(P'',Q) \in \ov{\gr}$ to $Q \in \ov{\img}$ is an isomorphism; and so, for each $Q \in \ov{\img}$, there exists a unique $P'' \in \mP^2$, such that $(P'',Q) \in \ov{\gr}$. Moreover, the coordinates of every point on $\gr$ satisfy $(\dag)$, so the coordinates of every point on $\ov{\gr}$ also satisfy $(\dag)$. Thus, the coordinates of $(P'',Q)$ satisfy $(\dag)$, and so all the $2 \times 2$ minors of $M$. Therefore, $P=P'=P''$, i.e. $(P,Q) \in \ov{\gr}$. We have shown that $\V \subseteq \ov{\gr}$. Hence, $\V = \ov{\gr}$.

In conclusion, the equations in $(\ref{3-sim1})$ and the $2 \times 2$ minors of $M$ describe $\ov{\gr}$ as a set. Furthermore, $M$ is a matrix of indeterminates, so it is a well known fact that the $2 \times 2$ minors of $M$ form a prime ideal. Similar to the last part of the proof of \cite[Theorem 3.6]{ha2}, we consider the following sequence of surjective ring homomorphisms:
\[ R[x_{ij}] \stackrel{\phi}{\rightarrow} R[w_it_j] \stackrel{\psi}{\rightarrow} R[w_iF_jt], \]
where $\phi$ sends $R$ to $R$, and sends $x_{ij}$ to $w_it_j$; and $\psi$ sends $R$ to $R$, and sends $w_it_j$ to $w_iF_jt$. Then from the proofs of equations $(\ref{3-break-1})$ and $(\ref{3-break-2})$, we further deduce that the $2 \times 2$ minors of $M$ form the kernel of $\phi$, and the images of equations in $(\ref{3-sim1})$ through $\phi$ form the kernel of $\psi$. Therefore, the $2 \times 2$ minors of $M$ and the equations in $(\ref{3-sim1})$ form the kernel of $\psi \circ \phi$, which is a prime ideal. Hence, the equations in $(\ref{3-sim1})$ and the $2 \times 2$ minors of $M$ form the defining ideal for $\ov{\gr}$ in $\mP^2 \times \mP^N$. The theorem is proved. 
\end{proof}

This gives rise to the following result.

\begin{thm} \label{3-simple-2}
Let $I = \oplus_{t \ge d} I_t$ be the defining ideal of $s = {d+1 \choose 2}$ points in $\mP^2$ which are in generic position. Then the defining equations for the Rees algebra $\R(I_{d+1})$ are the $2 \times 2$ minors of a $3 \times (d+2)$ matrix of linear forms. Moreover, $\R(I_{d+1})$ is Cohen-Macaulay, and has the same Betti numbers as that of the ideal of $2 \times 2$ minors of a generic $3 \times (d+2)$ matrix.
\end{thm}

\begin{proof}
The first statement of the theorem follows from Theorem \ref{3-simple-1} and the fact that the Rees algebra $\R(I_{d+1})$ is the bi-graded coordinate ring of $\ov{\Gamma_{d+1}}$ in $\mP^2 \times \mP^N$ (cf. \cite[Chapter 2]{ha-thesis}). For the second statement of the theorem, we observe that the defining ideal of $\ov{\gr}$ is the ideal of $2 \times 2$ minors of a matrix of linear forms of size $3 \times (d+2)$, so $\mbox{codim} \ov{\gr} \le 2(d+1)$. Furthermore, $\ov{\gr}$ is a surface in the product space $\mP^2 \times \mP^{2d+2}$ (after factoring out the linear forms in $(\ref{3-sim1})$), so its codimension is exactly $2(d+1)$. This implies that the defining ideal of $\ov{\gr}$ is perfect, and has the same Betti numbers as that of the ideal of $2 \times 2$ minors of a generic $3 \times (d+2)$ matrix. The result then follows.
\end{proof}

{\bf Remark:} The resolution of the ideal of minors of a generic matrix was computed by many authors (cf. \cite{la}, \cite{p-w}). One can apply their results to get the Betti numbers for $\R(I_{d+1})$. 

We lastly observe that it is not hard to extend the whole discussion to the case of reduced codimension 2 perfect ideals with linear presentation in a polynomial ring. More precisely, one can follow the same argument to obtain the following result.

\begin{thm} \label{3-simple-ext}
Suppose $I \subseteq R=\goth{k}[w_1,\ldots,w_n]$ is a reduced codimension 2 perfect ideal with a generic linear presentation. Let $I = \oplus_{t \ge d}I_t$ be its homogeneous decomposition. Then the Rees algebra $\R(I_{d+1})$ is Cohen-Macaulay, and its defining equations are the $2 \times 2$ minors of a $n \times (d+2)$ matrix of linear forms. Moreover, the Betti numbers of $\R(I_{d+1})$ are the same as those of the ideal of $2 \times 2$ minors of a generic $n \times (d+2)$ matrix.
\end{thm}

\begin{center}
{\bf Arbitrary number of points}
\end{center}

Our argument in this section inherits a great deal from that of \cite{gi-lo}. We refer the readers to Theorem 4.2 and Proposition 4.4 of \cite{gi-lo}. Suppose $s = {d+1 \choose 2} + k$ with $0 < k < d+1$. Then $\ix = \oplus_{t \ge d} I_t$, and $\sigma(\ix) = d+1$. The ideal generation conjecture is true in $\mP^2$ (cf. \cite{ggr} or \cite{g-m}), so we may take $\x$ general enough to have this conjecture satisfied. If we add the hypothesis that no $d+1$ points of $\x$ lie on a line, then the linear system $I_{d+1}$ is very ample on $\mP^2(\x)$ (\cite{d-g}). That is, the rational map from $\mP^2$ to $\mP^N$ given by a system of generators of the vector space $I_{d+1}$ gives an embedding of $\mP^2(\x)$. We shall now establish the defining equations for the Rees algebra $\R(I_{d+1})$ under these conditions.

The ideal generation conjecture states that $\ix$ is minimally generated by $d-k+1$ forms of degree $d$, and $h$ forms of degree $d+1$, where $h$ is either $0$ or $2k-d$, depending on whether $d \ge 2k$ or not. Moreover, by the Hilbert-Burch theorem, these generators can be seen as the $\rho+1$ maximal minors of a $\rho  \times (\rho+1)$ matrix $\bL$, where
\[ \rho = \left\{ \begin{array}{lc} k & \mbox{if } d \le 2k \\ d-k & \mbox{if } d \ge 2k. \end{array} \right. \]

In the case when $d < 2k$, the matrix $\bL$ is given by
\[ \bL = \left[ \begin{array}{cccccc} L_{1,1} & \ldots & L_{1,2k-d} & Q_{1,1} & \ldots & Q_{1,d-k+1} \\
\vdots &   & \vdots & \vdots &   & \vdots \\
L_{k,1} & \ldots & L_{k,2k-d} & Q_{k,1} & \ldots & Q_{k, d-k+1} \end{array} \right]. \]
where the $L_{i,j}$s are linear forms and the $Q_{i,j}$s are forms of degree 2. We denote by $F_j$ the minor obtained by deleting column $2k-d+j$ for $j=1,\ldots,d-k+1$, and by $G_l$ the minor obtained by deleting column $l$ for $l=1,\ldots,2k-d$. In this case $\ix = \ <F_1, \ldots, F_{d-k+1}, G_1,\ldots, G_{2k-d}>$, where $F_j \in I_d$ and $G_l \in I_{d+1}$ for all $j$ and $l$.

In the case when $d \ge 2k$, the matrix $\bL$ is given by
\[ \bL = \left[ \begin{array}{cccc} Q_{1,1} & Q_{1,2} & \ldots & Q_{1,d-k+1} \\
\vdots & \vdots & \ldots & \vdots \\
Q_{k,1} & Q_{k,2} & \ldots & Q_{k,d-k+1} \\
L_{1,1} & L_{1,2} & \ldots & L_{1,d-k+1} \\
\vdots & \vdots & \ldots & \vdots \\
L_{d-2k,1} & L_{d-2k,2} & \ldots & L_{d-2k,d-k+1} \end{array} \right]. \]
where, again, the $L_{i,j}$s are linear forms and the $Q_{i,j}$s are forms of degree 2. We denote by $F_j$ the minor obtained by deleting colomn $j$ for all $j=1,\ldots,d-k+1$. In this case, $\ix = \ <F_1,\ldots,F_{d-k+1}>$, where $F_j \in I_d$ for all $j$.

If $d < 2k$, a minimal system of generators for the vector space $I_{d+1}$ is given by $\{ w_iF_j, G_l | i=1,2,3; j=1,\ldots,d-k+1; l=1,\ldots,2k-d \}$. On the other hand, if $d \ge 2k$, a system of generators for the vector space $I_{d+1}$ is given by $\{ w_iF_j | i=1,2,3; j=1,\ldots,d-k+1 \}$, but this may not be minimal. In this case, the system $w_iF_j$'s gives $3(d-k+1)$ generators, while the vector space dimension of $I_{d+1}$ is $2d-k+1$, so there must be $d-2k$ linear dependence relations among those generators. Those relations can be found  by expanding the following zero determinants :
\[ 0 = \mbox{det } \left[ \begin{array}{c} \bL \\ L_{l,1} \ldots L_{l,d-k+1} \end{array} \right] = \left| \begin{array}{cccc} Q_{1,1} & Q_{1,2} & \ldots & Q_{1,d-k+1} \\ \vdots & \vdots & \ldots & \vdots \\
Q_{k,1} & Q_{k,2} & \ldots & Q_{k,d-k+1} \\
L_{1,1} & L_{1,2} & \ldots & L_{1,d-k+1} \\
\vdots & \vdots & \ldots & \vdots \\
L_{d-2k,1} & L_{d-2k,2} & \ldots & L_{d-2k,d-k+1} \\ L_{l,1} & L_{l,2} & \ldots & L_{l,d-k+1} \end{array} \right|. \]

Write $L_{l,j} = \sum_{i=1}^3 \lambda_{lji} w_i$ for all $l=1,\ldots,d-2k$, then
\begin{eqnarray} 
0 = \left| \begin{array}{c} \bL \\ L_{l,1} \ldots L_{l,d-k+1} \end{array} \right| & = & \sum_{j=1}^{d-k+1}L_{l,j}F_j = \sum_{j=1}^{d-k+1} \sum_{i=1}^{3} \lambda_{lji} w_iF_j. \label{3-temp}
\end{eqnarray}

Again, it was proved in \cite{gi-lo} that these $d-2k$ equations are linearly independent, so they are indeed all the dependence relations there are. 

We consider the rational map on $\mP^2$
\[ \varphi_{d+1} : \mP^2 \m \mP^N, \]
defined by sending each point $P=[\ov{w_1}:\ov{w_2}:\ov{w_3}] \not\in \x$ to $[\ov{w_i}F_j(P):G_l(P)]$ if $d < 2k$, or to $[\ov{w_i}F_j(P)]$ if $d \ge 2k$. Here, $(N+1)$ is the appropriate number of generators chosen for $I_{d+1}$. Again, let $\gr$ and $\img$ be the graph and the image of $\varphi_{d+1}$, and let $\ov{\gr}$ and $\ov{\img}$ be their closures in $\mP^2 \times \mP^N$ and $\mP^N$, respectively. If $d < 2k$, we use $[x_{ij}:y_l]$ to represent the homogeneous coordinates of $\mP^N$ such that $\varphi_{d+1}([\ov{w_1}:\ov{w_2}:\ov{w_3}]) = [\ov{x_{ij}}:\ov{y_l}]$, where $\ov{x_{ij}} = \ov{w_i}F_j(\ov{w_1}, \ov{w_2}, \ov{w_3})$ and $\ov{y_l} = G_l(\ov{w_1}, \ov{w_2}, \ov{w_3})$. If $d \ge 2k$, we use $[x_{ij}]$ to represent the homogeneous coordinates of $\mP^N$ such that $\varphi_{d+1}([\ov{w_1}:\ov{w_2}:\ov{w_3}]) = [\ov{x_{ij}}]$, where $\ov{x_{ij}} = \ov{w_i}F_j(\ov{w_1}, \ov{w_2}, \ov{w_3})$. 

When $d \ge 2k$, the linear dependence relations of the $w_iF_j$'s in $(\ref{3-temp})$ give rise to a collection of linear equations which are satisfied on $\gr$ as follows.
\begin{eqnarray}
\sum_{j=1}^{d-k+1} \sum_{i=1}^3 \lambda_{lji} x_{ij} & = & 0, \mbox{ for all } l=1,\ldots,d-2k. \label{3-linear}
\end{eqnarray}

Now, consider the matrix
\[ X = \left[ \begin{array}{cccc} w_1 & x_{11} & \ldots & x_{1, d-k+1} \\
w_2 & x_{21} & \ldots & x_{2, d-k+1} \\
w_3 & x_{31} & \ldots & x_{3, d-k+1} \end{array} \right]. \]
Clearly on $\gr$, the homogeneous coordinates of the points satisfy the $2 \times 2$ minors of $X$. To proceed, similar to what was done in \cite{gi-lo}, we separate the two cases.

{\bf Case 1:} $d < 2k$. In this case, for each $u=1,\ldots,k$, expanding the following zero determinant:
\[ \mbox{det } \bL_u = \left| \begin{array}{c} \bL \\ L_{u,1} \ldots L_{u,2k-d} ~ Q_{u,1} \ldots Q_{u,d-k+1} \end{array} \right|, \]
we obtain
\[ 0 = \sum_{l=1}^{2k-d} L_{u,l} G_l + \sum_{j=1}^{d-k+1} Q_{u,j}F_j. \]

Let $L_{u,l} = \sum_{i=1}^3 \lambda_{uli} w_i$ and $Q_{u,j} = \sum_{i,h=1}^3 \gamma_{uihj} w_iw_h$. Then
\[ 0 = \sum_{l=1}^{2k-d}(\sum_{i=1}^3 \lambda_{uli} w_i)G_l + \sum_{j=1}^{d-k+1}(\sum_{i,h=1}^3 \gamma_{uihj}w_iw_h)F_j. \]
Rewriting this as
\begin{eqnarray} 
0 & = & \sum_{i=1}^3 (\sum_{l=1}^{2k-d}\lambda_{uli}G_l)w_i + \sum_{i=1}^3(\sum_{j=1}^{d-k+1} \sum_{h=1}^3 \gamma_{uihj}w_hF_j)w_i. \label{3-eq1}
\end{eqnarray}

Also, for each $v=1,\ldots,d-k+1$, we get
\begin{eqnarray}
0 & = & F_v \bigg(\sum_{l=1}^{2k-d}(\sum_{i=1}^3 \lambda_{uli} w_i)G_l + \sum_{j=1}^{d-k+1}(\sum_{i,h=1}^3 \gamma_{uihj}w_iw_h)F_j\bigg) \nonumber \\
& = & \sum_{i=1}^3 (\sum_{l=1}^{2k-d}\lambda_{uli}G_l)w_iF_v + \sum_{i=1}^3(\sum_{j=1}^{d-k+1} \sum_{h=1}^3 \gamma_{uihj}w_hF_j)w_iF_v. \label{3-eq2}
\end{eqnarray}

The equations in (\ref{3-eq1}) and (\ref{3-eq2}) give a collection of bi-homogeneous equations that are satisfied by the coordinates of the points in $\gr$:
\[ \sum_{i=1}^3 (\sum_{l=1}^{2k-d}\lambda_{uli}y_l)w_i + \sum_{i=1}^3(\sum_{j=1}^{d-k+1}\sum_{h=1}^3\gamma_{uihj}x_{hj})w_i = 0, \]
for all $u=1,\ldots,k$, and 
\[ \sum_{i=1}^3 (\sum_{l=1}^{2k-d}\lambda_{uli}y_l)x_{iv} + \sum_{i=1}^3(\sum_{j=1}^{d-k+1}\sum_{h=1}^3\gamma_{uihj}x_{hj})x_{iv} = 0, \]
for all $u=1,\ldots,k$ and $v=1,\ldots,d-k+1$.  

Now, let $B=(b_{ui})_{1 \le u \le k, 1 \le i \le 3}$ be the matrix given by 
\[ b_{ui} = \sum_{l=1}^{2k-d} \lambda_{uli}y_l + \sum_{j=1}^{d-k+1} \sum_{h=1}^3 \gamma_{uihj}x_{hj}, \]
then the collection of the equations above can be rewritten as:
\[ \sum_{i=1}^3 b_{ui}w_i = 0, \mbox{ for any } u=1,\ldots,k, \]
and 
\[ \sum_{i=1}^3 b_{ui}x_{iv} = 0, \mbox{ for any } u=1,\ldots,k \mbox{ and } v=1,\ldots,d-k+1. \]
These equations are exactly all the entries of $B.X$ where $B$ and $X$ are the matrices as defined.

It follows also from \cite{gi-lo} that the coordinates of all the points of $\gr$ satisfy the $3 \times 3$ minors of $B$. 
We let ${\bf J}$ be the ideal in $\goth{k}[\under{w},\under{x},\under{y}]$ defined by the $3 \times 3$ minors of $B$, the $2 \times 2$ minors of $X$ and the entries of $B.X$.

{\bf Case 2:} $d \ge 2k$. Similar to what was done in the previous case, we expand the zero determinants:
\[ \mbox{det } \left[ \begin{array}{c} \bL \\ Q_{u,1} \ldots Q_{u,d-k+1} \end{array} \right], \]
for $u=1,\ldots,k$. We also let $B=(b_{ui})_{1 \le u \le k, 1 \le i \le 3}$ be the matrix given by
\[ b_{ui} = \sum_{j=1}^{d-k+1} \sum_{h=1}^3 \gamma_{uihj} x_{hj}, \]
where $Q_{u,j} = \sum_{i,h=1}^3 \gamma_{uihj} w_iw_h$. Again, let ${\bf J}$ be the ideal defined by the $3 \times 3$ minors of $B$, the $2 \times 2$ minors of $X$ and the entries of $B.X$. Similar to the previous case, we can prove that the coordinates of the points on $\gr$ satisfy the equations in ${\bf J}$.

We now establish a collection of equations which are similar to those in $(\dag)$ as follows:
\[ (\dag\dag) \left\{ \begin{array}{ccc} x_{2j}w_1 - x_{1j}w_2 & = & 0 \\ x_{3j}w_2-x_{2j}w_3 & = & 0 \\ 
x_{1j}w_3-x_{3j}w_1 & = & 0 \end{array} \right. \mbox{ for } j=1,2,\ldots,d-k+1. \]

From the proof of Proposition 4.4 of \cite{gi-lo}, one can deduce the following proposition:

\begin{pro} \label{3-pro2}
If $Q \in \mP^N$ whose homogeneous coordinates satisfy the $3 \times 3$ minors of $B$, the $2 \times 2$ minors of
\[ Y =  \left[ \begin{array}{ccc} x_{11} & \ldots & x_{1, d-k+1} \\
x_{21} & \ldots & x_{2, d-k+1} \\
x_{31} & \ldots & x_{3, d-k+1} \end{array} \right], \]
and the entries of $B.Y$, then there exists a unique point $P' \in \mP^2$ such that the homogeneous coordinates of $P'$ and $Q$ satisfy the equations in $(\dag\dag)$. 
\end{pro}

We obtain our first result when the number of points in $\x$ is arbitrary.

\begin{thm} \label{3-hard-1}
Let $\V$ be the subvariety of $\mP^2 \times \mP^N$ defined by ${\bf J}$ if $d<2k$, or defined by ${\bf J}$ and the equations in $(\ref{3-linear})$ if $d \ge 2k$. Then $\V = \ov{\gr}$ as sets.
\end{thm}

\begin{proof} The proof goes in the same manner of that of Theorem \ref{3-simple-1}. First, since all the points on $\gr$ satisfy the defining equations of $\V$, we have $\gr \subseteq \V$, whence $\ov{\gr} \subseteq \V$. 

To prove the reverse inclusion, suppose $(P,Q) \in \V$, where $P \in \mP^2$ and $Q \in \mP^N$. By Proposition \ref{3-pro2} and follow the same argument as that of \cite{gi-lo}, there exists a unique $P'=[\ov{w_1}:\ov{w_2}:\ov{w_3}]$ such that the coordinates of $P'$ and $Q$ satisfy the equations in $(\dag\dag)$, and $Q$ must have the form
\[ Q = [\ov{w_1}c_1:\ov{w_2}c_1:\ov{w_3}c_1:\ldots:\ov{w_1}c_{d-k+1}:\ov{w_2}c_{d-k+1}:\ov{w_3}c_{d-k+1}], \]
for some $c_1,\ldots,c_{d-k+1} \in \goth{k}$, if $d \ge 2k$, or 
\[ Q = [\ov{w_1}c_1:\ov{w_2}c_1:\ov{w_3}c_1:\ldots:\ov{w_1}c_{d-k+1}:\ov{w_2}c_{d-k+1}:\ov{w_3}c_{d-k+1}:\ov{y_1}:\ldots:\ov{y_{2k-d}}], \]
for some $c_1,\ldots,c_{d-k+1},\ov{y_1}, \ldots, \ov{y_{2k-d}} \in \goth{k}$, if $d < 2k$.

Now, substituting the coordinates of $Q$ into the entries of the product matrix $B.X$ and the equations in (\ref{3-linear}), we get 
\[ \bL (P') \left[ \begin{array}{c} c_1 \\ \vdots \\ c_{d-k+1} \end{array} \right] = \left[ \begin{array}{c} 0 \\ \vdots \\ 0 \end{array} \right], \]
if $d \ge 2k$, or 
\[ \bL (P') \left[ \begin{array}{c} c_1 \\ \vdots \\ c_{d-k+1} \\ y_1 \\ \vdots \\ y_{2k-d} \end{array} \right] = \left[ \begin{array}{c} 0 \\ \vdots \\ 0 \\ 0 \\ \vdots \\ 0 \end{array} \right], \]
if $d < 2k$. Thus, if $P' \not\in \x$, then 
\[ \left[ \begin{array}{c} c_1 \\ \vdots \\ c_{d-k+1} \end{array} \right] = \rho \left[ \begin{array}{c} F_1(P') \\ \vdots \\ F_{d-k+1}(P') \end{array} \right], \]
if $d \ge 2k$, or 
\[ \left[ \begin{array}{c} c_1 \\ \vdots \\ c_{d-k+1} \\ y_1 \\ \vdots \\ y_{2k-d} \end{array} \right] = \rho \left[ \begin{array}{c} F_1(P') \\ \vdots \\ F_{d-k+1}(P') \\ G_1(P') \\ \vdots \\ G_{2k-d}(P') \end{array} \right], \]
if $d < 2k$ (for some $\rho \in \goth{k}$); and if $P' \in \x$, then $Q$ belongs to the exceptional line corresponding the the blowup at $P'$. Therefore, $Q \in \ov{\img}$. 

We note further that the equations in $(\dag\dag)$ are among the equations of ${\bf J}$, so Proposition \ref{3-pro2} shows that for each $Q \in \ov{\img}$, there exists a unique $P'$ such that the coordinates of $P'$ and $Q$ satisfy the equations of ${\bf J}$ (when $d < 2k$), or of ${\bf J}$ and $(\ref{3-linear})$ (when $d \ge 2k$). Thus, $P = P'$. Now, with our assumption that there are no $d+1$ points of $\x$ lying on a line, the linear system $I_{d+1}$ is very ample, so the projection that sends $(P'',Q) \in \mP^2 \times \mP^N$ to $Q \in \mP^N$ is an isomorphism from $\ov{\gr}$ to $\ov{\img}$. Thus, for each $Q \in \ov{\img}$ there exists a unique $P'' \in \mP^2$ such that $(P'',Q) \in \ov{\gr}$. Moreover, since $\ov{\gr} \subseteq \V$, the coordinates of $P''$ and $Q$ satisfy the equations of ${\bf J}$ (when $d < 2k$), or of ${\bf J}$ and $(\ref{3-linear})$ (when $d \ge 2k$). Hence, $P = P' = P''$. In other words, $(P,Q) \in \ov{\gr}$. We just proved that $\V \subseteq \ov{\gr}$.

Hence, $\V = \gr$ as sets.
\end{proof}

This gives rise to the following result.

\begin{thm} \label{3-hard-2}
For a general choice of the points in $\x$, the defining equations of the Rees algebra $\R(I_{d+1})$ are the equations in ${\bf J}$ if $d < 2k$, or the equations in ${\bf J}$ together with the equations in $(\ref{3-linear})$ if $d \ge 2k$. Moreover, $\R(I_{d+1})$ is Cohen-Macaulay.
\end{thm}

\begin{proof} We first show that for a general choice of $\x$, the ideal ${\bf J}$ is prime and perfect. Similar to what was done in \cite{gi-lo}, we consider a new polynomial ring $R' = \goth{k}[\under{w},\under{x},\under{y},\under{z}]$ (in the case where $d \ge 2k, R'=\goth{k}[\under{w},\under{x},\under{z}]$), where $\under{z} = \{ z_{ui} \}_{1 \le u \le k, 1 \le i \le 3}$, and let $B' = (z_{ui})$. Then we can view ${\bf J}$ as the quotient ideal of ${\bf J}'$ in the ring $R'/(H_{ui})$, where ${\bf J}'$ is the ideal in $R'$ defined by the $3 \times 3$ minors of $B'$, the $2 \times 2$ minors of $X$ and the entries of $B'.X$, and
\[ H_{ui} = z_{ui} - b_{ui}, \mbox{ for all } u=1,\ldots,k \mbox{ and } i=1,2,3.\]

Since, ${\bf J}$ and ${\bf J}'$ are both bi-homogeneous, they are in particular also homogenous, so they define subvarieties of certain projective spaces. Let ${\bf W}$ be the subvariety of $\mP^{N+3}$ defined by ${\bf J}$, and let ${\bf W}'$ be the subvariety of $\mP^{N+3+3k}$ defined by ${\bf J}'$. Then ${\bf W}$ is obtained from ${\bf W}'$ by cutting ${\bf W}'$ with $3k$ hyperplanes $H_{ui}$. In other words, ${\bf W}$ is the intersection between ${\bf W}'$ and a $3k$-codimensional linear subspace of $\mP^{N+3+3k}$. By Huneke's theorem (\cite[Theorem 60]{hu}), we know that ${\bf W}'$ is an integral Cohen-Macaulay variety. 

Let $\Omega'$ be the grassmannian which parameterizes the linear subspaces of codimension $3k$ of $\mP^{N+3+3k}$. It follows from Bertini's theorem (cf. \cite{hart}, \cite{j}) that the subset ${\mathcal U}' \subseteq \Omega'$, such that for any $U' \in {\mathcal U}'$, $U' \cap {\bf W}'$ is again an integral Cohen-Macaulay variety, is non-empty and open. We also let $\Theta'$ be the grassmannian which parameterizes the linear subspaces of codimension $3k$ of $\mP^{N+3+3k}$ that lie inside the variety defined by the equations $w_1=w_2=w_3=0$. In \cite[Theorem 4.2]{gi-lo}, the authors actually showed that for a general choice of $\x$, the $3k$-codimensional linear subspace of $\mP^{N+3+3k}$ given by the hyperplanes $\{ H_{ui} | u=1,\ldots,k; i=1,2,3 \}$ is general in $\Theta'$. Moreover, consider the element $T' \in \Theta' \subseteq \Omega'$ given by $3k$ linear equations $z_{ui} = 0$, then $T' \cap {\bf W}'$ is the subvariety of $\mP^{N+3}$ given by the $2 \times 2$ minors of $X$, so it is an integral Cohen-Macaulay variety. In other words, $T' \in {\mathcal U}' \cap \Theta'$. Since $\Theta'$ is a closed subset of $\Omega'$, we deduce that ${\mathcal U}' \cap \Theta'$ is a non-empty open subset of $\Theta'$. All these facts, put together, imply that for a general choice of $\x$, the $3k$-codimensional linear subspace of $\mP^{N+3+3k}$ given by the hyperplanes $\{ H_{ui} | u=1,\ldots,k; i=1,2,3 \}$ is in ${\mathcal U}' \cap \Theta'$. In other words, for a general choice of $\x$, ${\bf W}$ is an integral Cohen-Macaulay subvariety of $\mP^{N+3}$, which, in turn, implies that ${\bf J}$ is a perfect prime ideal.

For $d < 2k$, this and Theorem \ref{3-hard-1} clearly imply that ${\bf J}$ is the defining ideal of $\R(I_{d+1})$. 

For $d \ge 2k$, let $l$ be the least integer such that $3l \ge d-2k$. Let $\Omega$ be the grassmannian which parameterizes the linear subspaces of codimension $3l$ of $\mP^{N+3}$, and let $\Theta$ be the grassmannian which parameterizes the linear subspaces of codimension $3l$ of $\mP^{N+3}$ that lie inside the variety defined by the equations $w_1 = w_2 = w_3 =0$. Again, it follows from Bertini's theorem that the subset ${\mathcal U} \subseteq \Omega$, such that for any $U \in {\mathcal U}$, $U \cap {\bf W}$ is an integral Cohen-Macaulay variety, is non-empty and open (for a general choice of $\x$, ${\bf W}$ is an integral Cohen-Macaulay variety). One can follow a similar argument as above, and consider the element $T$ in $\Theta$ given by $3l$ linear equations $\{ x_{ij} = 0 | i=1,2,3; j=d-k-l+2, \ldots, d-k+1 \}$, to show that ${\mathcal U} \cap \Theta$ is a non-empty open subset of $\Theta$. Moreover, from \cite{gi-lo}, it is easy to see that for a general choice of $\x$, the $3l$-codimensional linear subspace of $\mP^{N+3}$ given by the equations in (\ref{3-linear}) and $3l-(d-2k)$ other general hyperplanes is general in $\Theta$. Thus, for a general choice of $\x$, the $3l$-codimensional subspace of $\mP^{N+3}$ given by the equations in (\ref{3-linear}) and $3l-(d-2k)$ other general hyperplanes is in ${\mathcal U} \cap \Theta$, i.e. this $3l$-codimensional subspace of $\mP^{N+3}$ is general enough to intersect ${\bf W}$ at an integral Cohen-Macaulay variety. In particular, the hyperplanes given by the equations in (\ref{3-linear}) are general enough for a general choice of the points in $\x$. This and the fact in the previous paragraph show that for a general choice of $\x$, the hyperplanes $H_{ui}$ and the hyperplanes defined by the equations in (\ref{3-linear}) are general enough to intersect ${\bf W}'$ at an integral Cohen-Macaulay variety. Hence, for a general choice of $\x$, ${\bf J}$ together with the equations in $(\ref{3-linear})$ form a perfect prime ideal, i.e. ${\bf J}$ and the equations in $(\ref{3-linear})$ form the defining ideal of $\R(I_{d+1})$. 

The Cohen-Macaulayness of $\R(I_{d+1})$ follows from the perfection of its defining ideal. The theorem is proved.
\end{proof}

We finally observe that the same argument can be extended to a class of codimension 2 perfect ideals with presentation matrix in the form of $\bL$ (i.e. constituted by rows and columns of linear forms or quadratics). The generality of the points in $\x$ transforms to the genericity of the presentation matrix of the ideal. One can follow the same argument to obtain the following result.

\begin{thm} \label{3-hard-ext}
Suppose $I \subseteq R=\goth{k}[w_1,\ldots,w_n]$ is a generic codimension 2 perfect ideal, whose presentation matrix looks like that of $\bL$. Suppose also that $I = \oplus_{t \ge d} I_t$ is its homogeneous decomposition. Then the defining equations of the Rees algebra $\R(I_{d+1})$ are the $n \times n$ minors of a $k \times n$ matrix $B$ of linear forms, the $2 \times 2$ minors of an $n \times (d-k+2)$ matrix $X$ of linear forms, and the entries of $B.X$. Moreover, $\R(I_{d+1})$ is Cohen-Macaulay, and its defining ideal has the generic grade.
\end{thm}

\section{Asymptotic behaviour of the Rees algebras} \label{3-sec-asym}

If instead of a generic set of points in $\mP^2$, we start with an arbitrary set of points $\x$, then the presentation matrix of its defining ideal $\ix = \oplus_{t \ge \alpha} I_t$ no longer possesses a nice structure as it had in the previous section. It then becomes incredibly difficult to decide whether the Rees algebra $\R(I_t)$ is Cohen-Macaulay or to find its defining equations for a specific value of $t$. It is, however, possible to answer questions on the Cohen-Macaulayness or the degrees of the generators of $\R(I_t)$ as $t$ gets large. In this section, we address these questions. 

We begin by discussing a few properties of subschemes of a product scheme $\mP^n \times \mP^m$. For details on the definitions of product scheme $\mP^n \times \mP^m$, sheaves associated to bi-graded modules, and sheaf cohomology groups on $\mP^n \times \mP^m$, we refer the readers to \cite{vid}. Let $S = \goth{k}[x_0,\ldots,x_n,y_0,\ldots,y_m]$ be a polynomial ring over an algebraically closed field $\goth{k}$ of characteristic 0. Then $\mP^n \times \mP^m$, by definition, is the bi-Proj of $S$ with the natural bigradation on $S$ with respect to $(x_0,\ldots,x_n)$ and $(y_0,\ldots,y_m)$. Let $\goth{m} = (x_iy_j)_{i,j} \subseteq S$ be the bihomogeneous irrelevant ideal of $S$. We first recall the following known fact (cf. \cite{hy}, \cite{vid}).

\begin{pro} \label{3-1st-property}
Suppose $M$ is a bi-graded $S$-module, and $\sheaf{M}$ is the sheaf on $\mP^n \times \mP^m$ associated to $M$. Then, we have an exact sequence
\[ 0 \rightarrow H^0_{\goth{m}}(M) \rightarrow M \rightarrow \oplus_{a,b} H^0(\sheaf{M}(a,b)) \rightarrow H^1_{\goth{m}}(M) \rightarrow 0, \]
and isomorphisms
\[ \oplus_{a,b} H^i(\sheaf{M}(a,b)) \simeq H^{i+1}_{\goth{m}}(M) \ \forall \ i > 0. \]
\end{pro}

It is of our interest to study the Cohen-Macaulayness of the bi-graded coordinate ring of a variety in the product space $\mP^n \times \mP^m$. This is equivalent to the variety being {\it arithmetic Cohen-Macaulay} (a.CM). The definition of a.CM is given as follows.

\begin{mydef} Suppose $V \subseteq \mP^n \times \mP^m$ is a subscheme defined by the bihomogeneous ideal $I \subseteq S$. We say $V$ is {\it a.CM} if the bi-graded coordinate ring of $V$ is a Cohen-Macaulay ring, or equivalently, if $I$ is a perfect ideal in $S$.
\end{mydef}

On a product space $\mP^n \times \mP^m$ (for any $n$ and $m$), let $\pi_1 : \mP^n \times \mP^m \rightarrow \mP^n$ and $\pi_2 : \mP^n \times \mP^m \rightarrow \mP^m$ be the two projection maps. If $V$ is a subscheme of $\mP^n \times \mP^m$, when working on $V$, by abuse of notation, we also use $\pi_1$ and $\pi_2$ for the projection maps restricted on $V$. A necessary and a sufficient conditions for property a.CM are shown in the following theorem.

\begin{thm} \label{3-cm-condition}
Suppose $V \subseteq \mP^n \times \mP^m$ is a proper closed subscheme of dimension $d$ of $\mP^n \times \mP^m$. Then,
\renewcommand{\labelenumi}{(\arabic{enumi})}
\begin{enumerate}
\item If $V$ is a.CM, then $H^i(\mi_V(a,b)) = 0$ for all $a,b \in \z$ and $1 \le i \le d$, where $\mi_V$ is the ideal sheaf of $V$ in $\mP^n \times \mP^m$.
\item Suppose $d \not= n, m$, and $H^i(\mi_V(a,b)) = 0$ for all $a,b \in \z$ and $1 \le i \le d$. If in addition, $H^{d+1}(\mi_V(a,b)) = 0$ for all $a,b \ge 0$, and for every $j > 0$, 
\[ R^j {\pi_1}_{*}(\calo_V(p,q)) = 0 \ \forall p \in \z, q \ge 0, \]
and
\[ R^j {\pi_2}_{*}(\calo_V(p,q)) = 0 \ \forall q \in \z, p \ge 0, \]
\nopagebreak
then $V$ is a.CM.
\end{enumerate}
\end{thm}

\begin{proof} Similar criteria for varieties with negative $a^{*}$-invariants were given in \cite[Theorem 2.5]{hy}. We adopt his argument with a slight modification to prove our result. Let $\goth{n} = (x_0, \ldots, x_n, y_0, \ldots, y_m) \subseteq S$ be the maximal homogeneous ideal of $S$. Let $\goth{n}_1 = (x_0, \ldots, x_n)$ and $\goth{n}_2 = (y_0, \ldots, y_m)$ be the ideals in $S$ generated by the two sets of variables with respect to the standard bi-gradation of $S$. Then $\goth{n}_1 + \goth{n}_2 = \goth{n}$ and $\goth{n}_1 \cap \goth{n}_2 = \goth{m}$. Let $I_V$ and $S_V = S/I_V$ be the defining ideal and the coordinate ring of $V$, respectively. Then, $\dim S_V = d+2$. It is also not hard to see that the Cohen-Macaulayness of $S_V$ is equivalent to the condition that $H^i_{\goth{n}}(S_V) = 0$ for all $i=1, \ldots, d+1$, which is the same as the condition that $H^i_{\goth{n}}(I_V) = 0$ for all $i=1, \ldots, d+2$. 

{\bf (1)} Suppose that $V$ is a.CM. Equivalently, $H^i_{\goth{n}}(I_V) = 0$ for all $i=1, \ldots, d+2$. Consider the following Mayer-Vietoris sequence of local cohomology:
\[ \ldots \rightarrow H^i_{\goth{n}}(I_V) \rightarrow H^i_{\goth{n}_1}(I_V) \oplus H^i_{\goth{n}_2}(I_V) \rightarrow H^i_{\goth{m}}(I_V) \rightarrow H^{i+1}_{\goth{n}}(I_V) \rightarrow \ldots \]

The condition $H^i_{\goth{n}}(I_V) = 0$ for all $i=1, \ldots, d+2$, implies that the homomorphism
\[ H^i_{\goth{n}_1}(I_V) \oplus H^i_{\goth{n}_2}(I_V) \rightarrow H^i_{\goth{m}}(I_V) \]
is an isomorphism for $1 \le i \le d+1$ and injective for $i=d+2$. Localizing $H^i_{\goth{n}}(I_V)$ at the maximal ideals of $\oplus_{t \in \z}S_{(0,t)}$ and $\oplus_{t \in \z}S_{(t,0)}$, respectively, and making use of \cite[Lemma 1.1 and Lemma 2.3]{hy}, we have:
\[ H^i_{\goth{n}_1}(I_V) = H^i_{\goth{n}_2}(I_V) = 0 \ \forall i=1, \ldots, d+1. \]
This implies $H^i_{\goth{m}}(I_V) = 0$ for all $i=1,\ldots,d+1$. Together with Proposition \ref{3-1st-property}, it then follows that $H^i(\mi_V(a,b)) = 0$ for all $a,b \in \z$ and $i=1,\ldots,d$. 

{\bf (2)} Suppose now that $H^i(\mi_V(a,b)) = 0$ for all $a,b \in \z$ and $i=1,\ldots,d$, and $H^{d+1}(\mi_V(a,b)) = 0$ for all $a,b \ge 0$.  We observe the following. For $i=1$, this is to say that the homomorphism $S_{(a,b)} \rightarrow \Gamma(V, \calo_V(a,b))$ is surjective. In other words, the homomorphism 
\[ {S_V}_{(a,b)} \rightarrow \Gamma(V, \calo_V(a,b)) \]
is an isomorphism. Furthermore, the vanishing of $H^i(\mi_V(a,b))$ for all $a,b \in \z$, and all $i=1,\ldots,d$, is the same as having $H^i_{\goth{m}}(I_V) = 0$ for all $i=2,\ldots,d+1$. This is equivalent to having $H^i_{\goth{m}}(S_V) = 0$ for $i=1,\ldots,d$. Since $\goth{m} \not\subseteq I_V$, $H^0_{\goth{m}}(S_V)$ is clearly also 0. Lastly, the vanishing of $H^{d+1}(\mi_V(a,b))$ for all $a,b \ge 0$ implies that $H^d(\calo_V(a,b)) = 0$ for all $a,b \ge 0$, i.e. $[H^{d+1}_{\goth{m}}(S_V)]_{(a,b)} = 0$ for all $a,b \ge 0$.

Suppose, in addition, for every $j > 0$, $R^j {\pi_1}_{*}(\calo_V(p,q)) = 0 \ \forall p \in \z, q \ge 0$ and $R^j {\pi_2}_{*}(\calo_V(p,q)) = 0 \ \forall q \in \z, p \ge 0$. We need to show that $V$ is a.CM.

Let $T = \oplus_{t \in \z} {S_V}_{(t,0)}$ and $W = \mbox{Proj } T$, then $\pi_1 : V \rightarrow W$ is the canonical projection. For all $p \in \z$ and $q \ge 0$, the Leray spectral sequence 
\[ E^{i,j}_2 = H^i(W, R^j {\pi_1}_{*}(\calo_V(p,q))) \Rightarrow H^{i+j}(V, \calo_V(p,q)) \]
degenerates. Thus, the edge homomorphisms $H^i(W, {\pi_1}_{*}(\calo_V(p,q))) \rightarrow H^i(V, \calo_V(p,q))$ are just isomorphisms for all $p \in \z, q \ge 0$ and $i \ge 0$. 

Let $T_{<q>}$ ($q \ge 0$) be the $T$-module given by $T_{<q>} = \oplus_{t \in \z} {S_V}_{(t,q)}$, and let ${\mathcal T}_{<q>}$ be the sheaf associated to $T_{<q>}$ on $W$. It is easy to see that for $p \gg 0$, 
\[ \Gamma(W, {\mathcal T}_{<q>}(p)) = {S_V}_{(p,q)} = \Gamma(V, \calo_V(p,q)) = \Gamma(W, {\pi_1}_{*}(\calo_V(p,q))). \]
Also, ${\pi_1}_{*}(\calo_V(p,q)) \otimes \calo_W(p') = {\pi_1}_{*}(\calo_V(p+p',q))$. Thus, the canonical homomorphism ${\mathcal T}_{<q>}(p) \rightarrow {\pi_1}_{*}(\calo_V(p,q))$ is an isomorphism for all $p \in \z$ and $q \ge 0$. It now follows that the homomorphisms
\[ H^i(W, {\mathcal T}_{<q>}(p)) \rightarrow H^i(W, {\pi_1}_{*}(\calo_V(p,q))) \rightarrow H^i(V, \calo_V(p,q)), \]
similar to what was mentioned in \cite[Theorem 1.4]{hy}, are isomorphisms for all $p \in \z$, $q \ge 0$ and $i \ge 0$. Applying the five lemma on the diagram of \cite[Theorem 1.4]{hy}, it then implies that the homomorphism
\[ [H^i_{\goth{n}_1}(S_V)]_{(p,q)} \rightarrow [H^i_{\goth{m}}(S_V)]_{(p,q)} \]
is an isomorphism for all $p \in \z$, $q \ge 0$ and $i \ge 0$. By symmetry, the homomorphism
\[ [H^i_{\goth{n}_2}(S_V)]_{(p,q)} \rightarrow [H^i_{\goth{m}}(S_V)]_{(p,q)} \]
is also an isomorphism for all $q \in \z$, $p \ge 0$ and $i \ge 0$. Thus, for all $i =1, \ldots, d$, 
\[ [H^i_{\goth{n}_1}(S_V)]_{(p,q)} = 0 \ \forall p \in \z, q \ge 0 ~~ \mbox{ and } ~~ [H^i_{\goth{n}_2}(S_V)]_{(p,q)} = 0 \ \forall q \in \z, p \ge 0, \]
and
\[ [H^{d+1}_{\goth{n}_1}(S_V)]_{(p,q)} = 0 = [H^{d+1}_{\goth{n}_2}(S_V)]_{(p,q)} \ \forall p,q \ge 0. \]

Moreover, it is also easy to see that, for all $i > 0$,
\[ [H^i_{\goth{n}_1}(S_V)]_{(p,q)} = 0 \mbox{ if } q < 0 ~ ~ \mbox{ and } ~ ~ [H^i_{\goth{n}_2}(S_V)]_{(p,q)} = 0 \mbox{ if } p < 0. \]
Therefore, in the following Mayer-Vietoris sequence of local cohomology
\[ \ldots \rightarrow H^i_{\goth{n}}(S_V) \rightarrow H^i_{\goth{n}_1}(S_V) \oplus H^i_{\goth{n}_2}(S_V) \rightarrow H^i_{\goth{m}}(S_V) \rightarrow H^{i+1}_{\goth{n}}(S_V) \rightarrow \ldots \]
the homomorphisms
\[ H^i_{\goth{n}_1}(S_V) \oplus H^i_{\goth{n}_2}(S_V) \rightarrow H^i_{\goth{m}}(S_V) \]
are isomorphisms for $i=1, \ldots, d$, and injective for $i=d+1$. We get $H^i_{\goth{n}}(S_V) = 0$ for all $i=1,\ldots,d+1$. This is equivalent to $S_V$ being Cohen-Macaulay, i.e. $V$ being a.CM. The theorem is proved.
\end{proof}

Now, suppose $F \in S$ is a bihomogeneous polynomial. By abuse of notation, we denote by $(F)$ both the ideal generated by $F$ in $S$ and the subscheme of $\mP^n \times \mP^m$ defined by the equation $F = 0$.

\begin{pro} \label{3-2nd-property}
Suppose $V$ is a closed irreducible subscheme of $\mP^n \times \mP^m$ of dimension at least 2, and $L$ is a general linear form in the indeterminates $\{ y_j \ | \ j=0, \ldots, m \}$. Then $V$ is a.CM if and only if $V \cap (L)$, considered as a subscheme of $(L)$, is a.CM.
\end{pro}

\begin{proof} Let $I$ be the defining ideal of $V$. Then $I$ is a bihomogeneous ideal in $S$, so in particular, $I$ is a homogeneous ideal in $S$. Let $W$ then be the subscheme of $\mP^{n+m+1}$ defined by $I$. We observe that our discussion only involves the perfection of the defining ideal of $V$ (and $V \cap (L)$) which is the same as the defining ideal of $W$ (and $W \cap (L)$). Thus, the proposition would remain the same if instead of $V$ (and $V \cap (L)$) we look at $W$ (and $W \cap (L)$). The proposition now follows from \cite[Theorem 1.3.2]{mig}.
\end{proof}

From here onwards, we focus our attention back to the study of the asymptotic behaviour of the Rees algebras $\R(I_t)$ for the defining ideal of an arbitrary set of points in $\mP^2$. Again, suppose $\x = \{ P_1, \ldots, P_s \}$ is a set of $s$ distinct points in $\mP^2$. Let $\ix = \wp_1 \cap \ldots \cap \wp_s \subseteq R=\goth{k}[w_1,w_2,w_3]$ be the defining ideal of $\x$ ($\wp_i$ is the ideal of $P_i$), and suppose $\ix = \oplus_{t \ge \alpha} I_t$ is its homogeneous decomposition. 

\begin{thm} \label{3-cm}
Suppose $\x = \{ P_1, \ldots, P_s \}$ is an arbitrary set of $s$ points in $\mP^2$, and $\ix = \oplus_{t \ge \alpha} I_t \subseteq R = \goth{k}[w_1,w_2,w_3]$ is its defining ideal. Then, there exists an integer $d_0$ such that for all $t \ge d_0$, the Rees algebra $\R(I_t)$ of the ideal generated by $I_t$ is Cohen-Macaulay, and its defining ideal is generated by quadratics.
\end{thm}

\begin{proof} Let $\sigma = \sigma(\ix)$, the least integer at which the difference function of the Hilbert function of $\x$ equals 0, and take $d_0 = \max \{ 4, \sigma+1, s + 1 \}$ (note that $\sigma \le s$, so we in fact only need to take $d_0 = \max \{ 4, s+1 \}$). We shall prove that this value of $d_0$ satisfies the requirements of the theorem.

Suppose $t$ is an arbitrary integer which is bigger than or equal to $d_0$. We add ${t \choose 2} - s$ general smooth points to $\x$ to obtain a set of points $\tilde{\x}$ with the generic Hilbert function up to degree $t-2$, i.e.
\[ H_{\tilde{\x}} : 1 ~ 3 ~ 6 ~ \ldots ~ {t \choose 2} ~ {t \choose 2} \ \ldots. \]
We start by showing that the Rees algebra $\R(I_t)$ is Cohen-Macaulay using induction on the number $l = {t \choose 2} - s$ of points we add into $\x$ to get $\tilde{\x}$.

If $l = 0$, then $\x$ is a set of ${t \choose 2}$ points in $\mP^2$ with the generic Hilbert function. It follows from \cite{g-m} that the presentation matrix of $\ix$ has linear entries, $\ix$ is generated in degree $t-1$, and $\sigma(\ix) = t-1$. It now follows from Theorem \ref{3-simple-2} that the Rees algebra $\R(I_t)$ is Cohen-Macaulay. The assertion that the Rees algebra $\R(I_t)$ is Cohen-Macaulay is true for the base case.

Suppose now that our assertion is true for a set of points $\x' = \x \cup \{ P_{s+1} \}$, and we need to prove the assertion for the set of points $\x$. Let $\ix = \oplus_{t \ge \alpha} I_t$ and $\ix' = \oplus_{t \ge \alpha'} I_t'$ be the defining ideals of $\x$ and $\x'$ respectively. Since a general point $P_{s+1}$ imposes one independent condition at degree $t$, a system of generators for $I_t$ and for $I_t'$ may be given by $\{ F_0, \ldots, F_{r-1}, F_r \}$ and $\{ F_0, \ldots, F_{r-1} \}$, respectively. Consider the following rational maps
\[ \varphi : \mP^2 \m \mP^r \mbox{ and } \varphi' : \mP^2 \m \mP^{r-1} \]
given by $\varphi(P) = [F_0(P): \ldots: F_r(P)]$ and $\varphi'(P) = [F_0(P): \ldots: F_{r-1}(P)]$. Let $V$ and $V'$ be the closure of the graphs of these maps in $\mP^2 \times \mP^r$ and $\mP^2 \times \mP^{r-1}$, respectively. Clearly, the Rees algebras $\R(I_t)$ and $\R(I_t')$ are the bi-graded coordinate rings of $V$ and $V'$, respectively. By induction hypothesis, we know that $V'$ is a.CM. We need to show that $V$ is also a.CM.

Let $[y_0:\ldots:y_r]$ represent the homogeneous coordinates of $\mP^r$. Let $S = \goth{k}[w_1,w_2,w_3,\linebreak y_0, \ldots, y_r]$ and $S_V$ be the bi-graded coordinate rings of $\mP^2 \times \mP^r$ and $V$, respectively. Let $H = V \cap (y_r)$. By Proposition \ref{3-2nd-property}, we only need to show that $H$, considered as a subscheme of $\mP^2 \times \mP^{r-1}$, is a.CM. 

Clearly, $\pi_1(V) = \pi_1(V') = \mP^2$. Let $\ov{V} = \pi_2(V)$ and $\ov{V'} = \pi_2(V')$. By the construction of $\tilde{\x}$, we know that $\sigma(\ix) \le \sigma(I_{\tilde{\x}}) = t-1$ and $\sigma(I_{\x'}) \le \sigma(I_{\tilde{\x}}) = t-1$. Therefore, the linear systems $I_t$ and $I_t'$ are very ample (\cite{d-g}). Thus, $\pi_2$ is an isomorphism from $V$ onto $\ov{V}$ and from $V'$ onto $\ov{V'}$. It is easy to see that since the coordinate $y_r$ of $\mP^r$ is chosen generally, $H$ meets each exceptional curve of $V$ no more than once. Thus, $\pi_1^{-1}(P)$, restricted to $H$, is at most one point, for every $P \in \mP^2$. Since $\pi_2$ is an isomorphism on $V$, it is also an isomorphism on $H$, whence $\pi_2^{-1}(Q)$, restricted to $H$, is also at most one point, for every $Q \in \mP^{r-1}$. Therefore, by \cite[Corollary III.11.2]{hart}, for every $j > 0$, one gets
\[ R^j {\pi_1}_{*}(\calo_H(p,q)) = 0 ~ \mbox{ and } ~ R^j {\pi_2}_{*}(\calo_H(p,q)) = 0 ~ ~ \forall p,q \in \z. \]
By Theorem \ref{3-cm-condition}, to show that $H$ is a.CM, it is now enough to show $H^1(\mi_H(a,b)) = 0$ for all $a,b \in \z$, and $H^2(\mi_H(a,b)) = 0$ for all $a,b \ge 0$, where $\mi_H$ is the ideal sheaf of $H$ in $\mP^2 \times \mP^{r-1}$. 

It is easy to see that since $V'$ is the projection of $V$ on $(y_r)$ centered at the point $P_{s+1} \times \varphi(P_{s+1})$, $H \subseteq V'$. Consider the following exact sequence
\[ 0 \rightarrow \mi_{V'} \rightarrow \mi_H \rightarrow \mi_{H, V'} \rightarrow 0, \]
where $\mi_{H, V'}$ is the ideal sheaf of $H$ considered as a subscheme of $V'$. Since $V'$ is a.CM, taking the cohomology groups, we get
\[ H^1(\mi_H(a,b)) = H^1(\mi_{H, V'}(a,b)) ~~ \mbox{ and } ~~ H^2(\mi_H(a,b)) \hookrightarrow H^2(\mi_{H, V'}(a,b)). \]
Thus, it suffices to show that
\[ H^1(\mi_{H,V'}(a,b)) = 0 \ \forall a,b \in \z ~~ \mbox{ and } ~~ H^2(\mi_{H, V'}(a,b)) = 0 \ \forall a,b \ge 0. \]

Let $E_0$ be the pull back to $V'$ of the class of a general line in $\mP^2$, and let $E_1, \ldots, E_{s+1}$ be the classes of the exceptional divisors corresponding to the blowup at the points $P_1, \ldots, P_{s+1}$, respectively. Let $\ov{E_0}, \ldots, \ov{E_{s+1}}$ be the images of $E_0, \ldots, E_{s+1}$ through $\pi_2$, respectively, then they generate the Picard group of $\ov{V'}$. $\ov{V} \cap (y_r)$ is a hyperplane section of $\ov{V}$ and since the coordinates in $\mP^r$ are chosen generally, we may assume that $\ov{V} \cap (y_r)$ belongs to the divisor class $| t\ov{E_0} - \ov{E_1} - \ldots - \ov{E_s} |$. Thus, $H \in | tE_0 - E_1 - \ldots - E_s |$. We have
\begin{eqnarray*}
\lefteqn{\mi_{H, V'}(a,b) = \calo_{V'}(-H)(a,b)} \\
& = & \calo_{V'}(-H) \otimes \pi_1^{*} (\calo_{\mP^2}(a)) \otimes \pi_2^{*} (\calo_{\ov{V'}}(b)) \\
& = & \calo_{V'}(-H) \otimes \pi_1^{*} (\calo_{\mP^2}(a)) \otimes \pi_2^{*} (\calo_{\ov{V'}}(bt\ov{E_0} - b\ov{E_1} - \ldots - b\ov{E_s} - b\ov{E_{s+1}})) \\
& = & \calo_{V'}(-H) \otimes \calo_{V'}(aE_0) \otimes \calo_{V'}(btE_0 - bE_1 - \ldots - bE_s - bE_{s+1}) \\ 
& = & \calo_{V'}( (a+(b-1)t)E_0 - (b-1)E_1 - \ldots - (b-1)E_s - bE_{s+1})
\end{eqnarray*}

Let $D_{a,b} = (a+(b-1)t)E_0 - (b-1)E_1 - \ldots - (b-1)E_s - bE_{s+1}$ on $V'$. We first prove $H^2(\calo_{V'}(D_{a,b})) = 0$ for all $a,b \ge 0$. We shall use double induction on $b$ and $a$. 

For $b = 0$ and $0 \le a \le t$, we first have $D_{a,0} = (a-t)E_0 + E_1 + \ldots + E_s$. The canonical divisor on $V'$ is $K_{V'} = -3E_0 + E_1 + \ldots + E_{s+1}$. Let $H' = tE_0 - E_1 - \ldots - E_{s+1}$. From \cite{d-g}, $H'$ is very ample on $V'$. We also have, $K_{V'}.D_{a,0} = 3(t-a) - s > -3t + s + 1 = K_{V'}.H'$. Thus, $H^2(\calo_{V'}(D_{a,0})) = 0$ (cf. \cite[Lemma V.1.7]{hart}). Suppose now that $H^2(\calo_{V'}(D_{a,0})) = 0$ is true for $a \ge t$, we shall show that $H^2(\calo_{V'}(D_{a+1, 0})) = 0$. Indeed, consider the exact sequence
\[ 0 \rightarrow \calo_{V'}(D_{a,0}) \rightarrow \calo_{V'}(D_{a+1,0}) \rightarrow \calo_{E_0}(D_{a+1,0}) \rightarrow 0. \]
Since $\deg \calo_{E_0}(D_{a+1,0}) = a+1-t > 0$, and since $E_0$ is a rational curve, we have 
\[ H^1(\calo_{E_0}(D_{a+1,0})) = 0 \mbox{ and } H^2(\calo_{E_0}(D_{a+1,0})) = 0. \]
Thus, $H^2(\calo_{V'}(D_{a+1,0})) = H^2(\calo_{V'}(D_{a,0})) = 0$.

For $b=1$, we first have $D_{0,1} =  - E_{s+1}$. We also have $K_{V'}.D_{0,1} = 1 > -3t + s + 1 = K_{V'}.H'$. Thus, $H^2(\calo_{V'}(D_{0,1})) = 0$ (cf. \cite[Lemma V.1.7]{hart}). Suppose now that $H^2(\calo_{V'}(D_{a,1})) = 0$ is true for $a \ge 0$, we shall show that $H^2(\calo_{V'}(D_{a+1,1})) = 0$. Indeed, as before, consider the exact sequence
\[ 0 \rightarrow \calo_{V'}(D_{a,1}) \rightarrow \calo_{V'}(D_{a+1,1}) \rightarrow \calo_{E_0}(D_{a+1,1}) \rightarrow 0. \]
Since $\deg \calo_{E_0}(D_{a+1,1}) = a+1 > 0$, and since $E_0$ is a rational curve, we have 
\[ H^1(\calo_{E_0}(D_{a+1,1})) = H^2(\calo_{E_0}(D_{a+1,1})) = 0. \]
Thus, $H^2(\calo_{V'}(D_{a+1,1})) = H^2(\calo_{V'}(D_{a,1})) = 0$.

Now, suppose that $H^2(\calo_{V'}(D_{a,b})) = 0$ is true for any $a \ge 0$ and some $b \ge 1$, we shall show that $H^2(\calo_{V'}(D_{a,b+1})) = 0$. Indeed, consider the exact sequence
\[ 0 \rightarrow \calo_{V'}(D_{a,b}) \rightarrow \calo_{V'}(D_{a,b+1})  \rightarrow \calo_{H'}(D_{a,b+1}) \rightarrow 0. \]
Since $\deg \calo_{H'}(D_{a,b+1})  = (a+bt)t - b(s+1) - 1 > 2{t-1 \choose 2}+1 = 2g(H') + 1$, we have $H^1(\calo_{H'}(D_{a,b+1})) = H^2(\calo_{H'}(D_{a,b+1})) = 0$. Thus, 
\[ H^2(\calo_{V'}(D_{a,b+1})) = H^2(\calo_{V'}(D_{a,b})) = 0. \]

Hence, $H^2(\calo_{V'}(D_{a,b})) = 0$ for all $a,b \ge 0$.

It remains now to prove that $H^1(\calo_{V'}(D_{a,b})) = 0$ for all $a,b \in \z$. Similar to what we did above, we can rewrite 
\[ \calo_{V'}(D_{a,b}) = \calo_{V'}(-E_{s+1})(a,b-1) = \mi_{E_{s+1}, V'}(a,b-1), \]
where $\mi_{E_{s+1}, V'}$ is the ideal sheaf of $E_{s+1}$ considered as a subscheme of $V'$. It was proved in \cite[Proposition 2.1]{gi2} that $H^1(\calo_{\ov{V'}}((n-1)t\ov{E_0} - (n-1)\ov{E_1} - \ldots - (n-1)\ov{E_s} - n\ov{E_{s+1}}) = 0$ for all $n \in \z$. That is, $H^1(\mi_{\ov{E_{s+1}}, \ov{V'}}(n-1)) = 0$ for all $n \in \z$, where $\mi_{\ov{E_{s+1}}, \ov{V'}}$ is the ideal sheaf of $\ov{E_{s+1}}$ considered as a subscheme of $\ov{V'}$. Therefore, $\ov{E_{s+1}}$, considered as a subscheme of $\ov{V'}$, is projective CM (see \cite{gi2} for definition of projective CM). By considering the exact sequence
\[ 0 \rightarrow \mi_{\ov{V'}} \rightarrow \mi_{\ov{E_{s+1}}} \rightarrow \mi_{\ov{E_{s+1}}, \ov{V'}} \rightarrow 0, \]
where $\mi_{\ov{V'}}$ and $\mi_{\ov{E_{s+1}}}$ are the ideal sheaves of $\ov{V'}$ and $\ov{E_{s+1}}$ in $\mP^{r-1}$, and from the fact that $\ov{V'}$ is projective CM (\cite{gi2}), we deduce that in $\mP^{r-1}$, the homogeneous coordinate ring of $\ov{E_{s+1}}$ is CM. It is also clear that $E_{s+1} = P_{s+1} \times \ov{E_{s+1}}$. Thus, the coordinate ring of $E_{s+1}$ in $\mP^2 \times \mP^{r-1}$ is CM, i.e. $E_{s+1}$ is a.CM in $\mP^2 \times \mP^{r-1}$. Hence, by Theorem \ref{3-cm-condition}, $H^1(\mi_{E_{s+1}}(a,b-1)) = 0$ for all $a,b \in \z$, where $\mi_{E_{s+1}}$ is the ideal sheaf of $E_{s+1}$ in $\mP^2 \times \mP^{r-1}$. By considering the exact sequence
\[ 0 \rightarrow \mi_{V'} \rightarrow \mi_{E_{s+1}} \rightarrow \mi_{E_{s+1}, V'} \rightarrow 0, \]
and from the hypothesis that $V'$ is a.CM, we conclude that $H^1(\mi_{E_{s+1}, V'}(a,b-1)) = 0$ for all $a,b \in \z$. Hence,
\[ H^1(\calo_{V'}(D_{a,b})) = 0 \mbox{ for all } a,b \in \z. \]

The assertion is proved, i.e. $\R(I_t)$ is Cohen-Macaulay for all $t \ge d_0$.

We now proceed to prove that the defining ideal of $\R(I_t)$ is generated by quadratics. Again, we use induction on $l$, the number of general smooth points being added to $\x$ to get $\tilde{\x}$. If $l=0$, then it follows from Theorem \ref{3-simple-2} that the defining ideal of $\R(I_t)$ is generated by quadratics. Thus, the base case of the assertion that the defining ideal of the Rees algebra $\R(I_t)$ is generated by quadratics is proved. 

Suppose now that this assertion is true for a set of points $\x' = \x \cup \{ P_{s+1} \}$, and we need to show it for the set of points $\x$. Let $V, V'$ and $H$ be defined as before. By induction hypothesis, $V'$ is generated by quadratics. We need to show that $V$ is also generated by quadratics. Since $V$ is a.CM, this amounts to showing that $H$ is generated by quadratics. Moreover, since $H \subseteq V'$ and $V'$ is generated by quadratics, we only need to show that $I_{H, V'}$ (the ideal of $H$ in the coordinate ring of $V'$) is generated in its quadratic degree. It suffices to show that the following multiplication maps are onto:
\begin{eqnarray}
& & H^0(\mi_{H, V'}(2,0)) \otimes H^0(\calo_{V'}(1,0)) \smap H^0(\mi_{H, V'}(3,0)) \label{3-30} \\
& & H^0(\mi_{H, V'}(0,2)) \otimes H^0(\calo_{V'}(0,1)) \smap H^0(\mi_{H, V'}(0,3)) \label{3-03} \\
& & H^0(\mi_{H, V'}(1,1)) \otimes H^0(\calo_{V'}(1,0)) \smap H^0(\mi_{H, V'}(2,1)) \label{3-21} \\
& & H^0(\mi_{H, V'}(0,2)) \otimes H^0(\calo_{V'}(1,0)) \smap H^0(\mi_{H, V'}(1,2)), \label{3-12} 
\end{eqnarray}
where $\mi_{H, V'}$ is the ideal sheaf of $H$ in $V'$.

To prove (\ref{3-30}), for each integer $a$, as it was done above, we rewrite 
\[ \mi_{H, V'}(a,0) = \calo_{V'}(-H)(a,0) = \calo_{V'}((a-t)E_0+E_1+\ldots+E_s). \]
Again, the canonical divisor on $V'$ is $K_{V'} = -3E_0 + E_1+\ldots+E_s+E_{s+1}$. Let 
\[ D_a = (a-t)E_0+E_1+\ldots+E_s ~ \mbox{ and } ~ H' = tE_0-E_1-\ldots-E_s-E_{s+1}. \]
As before, from \cite{d-g}, $H'$ is every ample on $V'$. For $a=2,3$, we also have 
\[ (K_{V'} - D_a).H' = (t-a-3)t + 1 > -3t + s + 1 = K_{V'}.H'. \]
Therefore, it follows from \cite[Lemma V.1.7]{hart} that, for $a=2 \mbox{ and } 3$, 
\[ H^0(\calo_{V'}(D_a)) = H^2(\calo_{V'}(K_{V'}-D_a)) = 0. \] 
Hence, (\ref{3-30}) follows vacuously.

To prove (\ref{3-03}), for each integer $b$, as before, we rewrite
\begin{eqnarray*}
\mi_{H, V'}(0,b) & = & \calo_{V'}((b-1)tE_0 - (b-1)E_1 - \ldots - (b-1)E_s - bE_{s+1}) \\
& = & \mi_{E_{s+1}, V'}(0,b-1). 
\end{eqnarray*}
Thus, (\ref{3-03}) is equivalent to 
\[ H^0(\mi_{E_{s+1}, V'}(0,1)) \otimes H^0(\calo_{V'}(0,1)) \smap H^0(\mi_{E_{s+1}, V'}(0,2)). \]
This is indeed true, since $E_{s+1}$ is a line in $V'$, and hence, is obviously generated in degree 1.

Similarly, it can be shown that (\ref{3-21}) and (\ref{3-12}) are equivalent to
\[ H^0(\mi_{E_{s+1}, V'}(1,0)) \otimes H^0(\calo_{V'}(1,0)) \smap H^0(\mi_{E_{s+1}, V'}(2,0)), \]
and 
\[ H^0(\mi_{E_{s+1}, V'}(0,1)) \otimes H^0(\calo_{V'}(1,0)) \smap H^0(\mi_{E_{s+1}, V'}(1,1)), \]
respectively. Those are again true since $E_{s+1}$ is a line on $V'$, and so, generated in degree 1. The theorem is proved.
\end{proof}

{\bf Remark:} Most of our arguments do not use the fact that $\x$ is a set of points. Thus, if we couple with \cite[Theorem 1.3]{m-u}, the Cohen-Macaulayness of Theorem \ref{3-cm} is actually true for any subscheme of $\mP^2$ whose defining ideal satisfies condition $G_3$. An example of such subscheme is any locally complete intersection subscheme of $\mP^2$. 

Inspired by the notion of having property $N_p$, introduced by Green, together with many supportive evidences, we conjecture the following.

{\bf Conjecture:} Suppose $\x$ is an arbitrary set of points in $\mP^2$, and $\ix = \oplus_{t \ge \alpha} I_t \subseteq R = \goth{k}[w_1,w_2,w_3]$ its defining ideal. Then, for any non-negative integer $p$, there exists an integer $d_p$ such that for any $t \ge d_p$, the Rees algebra $\R(I_t)$ is Cohen-Macaulay, generated by quadratics, and the first $p$ steps in its minimal free resolution are linear, i.e. the first $p$ matrices are of linear forms.

Theorem \ref{3-cm} proves this conjecture for $p = 0$. 

\begin{small}
{\it Acknowledgement}. {\sf This paper is part of the author's PhD thesis. I would like to thank my research advisor A.V. Geramita for his inspiration and guidance. I would also like to thank Prof. A. Iarrobino, Prof. M. Johnson, Dr. J. Chipalkatti and Mr. M. Mustata for their useful suggestions during the preparation of this paper.}
\end{small}

\end{document}